\documentclass[graybox]{svmult}

\usepackage{dsfont}
\usepackage{mathptmx}       
\usepackage{helvet}         
\usepackage{courier}        
\usepackage{type1cm}        
\usepackage{makeidx}         
\usepackage{graphicx}        
\usepackage{multicol}        
\usepackage[bottom]{footmisc}

\makeindex             

\begin{document}

\title*{Global weak solutions of PDEs for compressible media:
         A compactness criterion to cover  new physical situations}
\titlerunning{A compactness criterion to cover new physical situations}
\author{D. Bresch, P.--E. Jabin}
\institute{Didier Bresch \at LAMA CNRS UMR 5127, University of Savoie Mont-Blanc, Bat. Le Chablais, Campus scientifique, 73376 Le Bourget du Lac, France. D.~Bresch is partially supported by the ANR- 13-BS01-0003-01 project DYFICOLTI., \email{didier.bresch@univ-smb.fr}
\and Pierre--Emmanuel Jabin \at CSCAMM and Dept. of Mathematics, University of Maryland,
College Park, MD 20742, USA. P.--E.~Jabin is partially supported by NSF Grant 1312142 and by NSF Grant RNMS (Ki-Net) 1107444. \email{pjabin@cscamm.umd.edu}}
%
%
\maketitle

\abstract*{This short paper is an introduction of  the memoir  recently written by the two
 authors (see [D.Bresch., P.--E. Jabin, arXiv:1507.04629, (2015)]) which concerns the resolution of two longstanding problems:
 Global existence of weak solutions for compressible Navier--Stokes equations 
 with  {\it  thermodynamically unstable pressure} and 
 with {\it anisotropic stress tensor}.  
 We focus here on a Stokes-like system which can for instance 
 model flows in a compressible tissue in biology or in a compressible porous
 media  in petroleum  engineering.
   This allows us to explain, on a simpler but still relevant and important system,  the tools recently 
 introduced by the authors and to discuss the important results that have been 
 obtained on  the compressible Navier--Stokes equations.
    It is finally a real pleasure to dedicate this paper to G. {\sc M\'etivier} for his 65's Birthday.}

\abstract{This short paper is an introduction of  the memoir  recently written by the two
 authors (see [D.Bresch., P.--E. Jabin, arXiv:1507.04629, (2015)]) which concerns the resolution of two longstanding problems:
 Global existence of weak solutions for compressible Navier--Stokes equations 
 with  {\it  thermodynamically unstable pressure} and 
 with {\it anisotropic stress tensor}.  
 We focus here on a Stokes-like system which can  for instance 
 model flows in a compressible tissue in biology or in a compressible porous
 media  in petroleum  engineering.
   This allows to explain, on a simpler but still relevant and important system,  the tools recently 
 introduced by the authors and to discuss the important results that have been 
 obtained on  the compressible Navier--Stokes equations. It is finally a real pleasure to dedicate
 this paper to G. {\sc M\'etivier} for his 65's Birthday.}

\section{Introduction}
    We consider in this paper a model which has been developed for
 flows in a compressible tissue in biology (see \cite{BrCoGrRiSa}, \cite{DoTr}) 
 or in compressible porous media  in petroleum engineering (see \cite{GaTo}).  
    The most simple system involves a density 
 $\rho$ that is transported, 
\[ 
\partial_t \rho + {\rm div} (\rho u) = 0,
\]
by a  velocity field $u$ described by a Stokes-like equation
\[
- \mu \Delta u + \alpha u +  \nabla P(\rho)= S, 
\]
with $\mu, \alpha>0$.  
    For simplicity we consider periodic boundary conditions, namely both equations are posed for
$x\in \Omega = {\mathds T}^d$. This is also the reason for the damping term $\alpha u$ to control $u$ without imposing any additional condition on $S$. The corresponding
PDE is usually named Brinkman equation. It accounts for flow through medium where the grains are porous themselves. 

   In this short paper,  we explain how to consider non-monotone pressure laws $P$ for
this system (complex pressure laws (attractive and repulsive)) to obtain the existence of global  weak-solutions. Note that in particular biological systems frequently exhibit preferred ranges of  densities for instance attractive interactions for low densities and repulsive at higher ones.
    
To get such global existence of weak solutions result, the two authors have
recently revisited (see \cite{BrJa})  the classical compactness theory on the density by obtaining precise quantitative regularity estimates: This requires a more precise analysis of the structure of the equations combined to a novel approach to the compactness of the continuity equation
 (by introducing appropriate weights).
We quote at the end of the article some of the precise results obtained in \cite{BrJa} on the compressible Navier-Stokes systems but we of course refer the reader to \cite{BrJa} for all the details and possible extensions for instance including temperature conductivity dependency.
\section{Equations and main result}
As mentioned above, we work on the torus ${\mathds T}^d$. This is only for simplicity in order to avoid discussing boundary conditions or the behavior at infinity.

\subsection{Statements of the result}
   We present in this section our main existence result concerning System (\ref{Biononmono}). As usual for global existence of weak solutions to nonlinear PDEs, one has to prove stability estimates for sequences of approximate solutions and construct such approximate sequences. The main contribution in this paper and the major part of the proofs concern the stability procedure and more precisely the compactness of the density.  We refer to \cite{BrJa} for details and the way to construct the approximate solutions sequence.
As per the introduction, we consider the following system
\begin{equation}\label{Biononmono}
\left\{
\begin{array}{rl}
& \partial_t \rho + {\rm div} (\rho u) =0, \\
& - \mu \Delta u  + \alpha u +  \nabla P(\rho) = S,\\
\end{array}
\right. 
\end{equation} 
with $\mu, \alpha >0$, a pressure law $P$ which is continuous on $[0,+\infty)$, $P$ locally Lipschitz on $(0,+\infty)$ with $P(0)=0$ such that there exists $C>0$ with
\begin{equation}\label{gammacontrolBrinkman}
C^{-1} \rho^\gamma - C  \le P(\rho) \le C \rho^\gamma + C,
\end{equation}
and for all $s\ge 0$
\begin{equation}
|P'(s)|\leq \bar P s^{\gamma-1}.
\label{nonmonotonehypBrinkman}
\end{equation}
 System {\rm (\ref{Biononmono})}  is completed with the initial boundary condition:
\begin{equation}\label{init}
\rho\vert_{t=0} = \rho_0.
\end{equation}
   One then has global existence of a weak solution.
\begin{theorem} \label{MainResultPressureLaw} 
Assume that $S\in L^2(0,T;\;H^{-1}({\mathds T}^d))$ and the periodic initial data $\rho_0$ satisfies  the bound
\[
\rho_0\ge 0, \qquad  0<  M_0 = \int_{{\mathds T}^d} \rho_0 < + \infty, \qquad
E_0= \int_{{\mathds T}^d}  \rho_0 e(\rho_0) \,dx <+\infty, 
\]
where $e(\rho) = \int_{\rho^\star}^\rho P(s)/s^2 ds$ with $\rho^\star$ a 
constant reference density.
   Let the pressure law $P$  satisfy
{\rm (\ref{gammacontrolBrinkman})} and {\rm (\ref{nonmonotonehypBrinkman})} with 
$\gamma > 1$.  
  Then there exists a global  weak solution $(\rho,u)$ of the compressible system
with positive density  in the sense that it satisfies the energy estimate  {\rm (\ref{ineg})}, 
the estimates
$$ \rho \in  \bigl(L^\infty(0,T;L^\gamma({\mathds T}^d)) \cap \,  L^{2\gamma}((0,T)\times {\mathds T}^d)\bigr)
                   \cap \,  C([0,T]; L^\gamma({\mathds T}^d) \hbox{ weak}),$$
$$ u \in L^2(0,T;H^1({\mathds T}^d))$$
and Equations    {\rm (\ref{Biononmono})} and {\rm (\ref{init})} 
respectively  in ${\cal D}'((0,T)\times {\mathds T}^d)$ and in ${\cal D}'({\mathds T}^d)$.
\end{theorem}
 
 \bigskip

\begin{remark}
   From the bounds provided by the theorem, it is straightforward to check that $\rho\,u\in L^1((0,T)\times{\mathds T}^d)$. Similarly from (\ref{gammacontrolBrinkman}) and the theorem, $P(\rho)\in L^1((0,T)\times{\mathds T}^d)$. Therefore all the terms in Equations  {\rm (\ref{Biononmono})} and {\rm (\ref{init})} make sense in ${\cal D}'((0,T)\times {\mathds T}^d)$. Note that since all the terms in the second equation in {\rm (\ref{Biononmono})} are in $L^1_{loc}$ in time, this equation could even be posed for $a.e.\;t$.
   A weak formulation of   {\rm (\ref{Biononmono})}  may also be written as usually for global weak solutions "\`a la J. {\sc Leray}".
\end{remark} 

 \begin{remark}
Let us note that we do not try to optimize the regularity of $S$ which could be far less smooth.
The objective of this short note being to be an introduction to \cite{BrJa}
focusing on the new compactness criterion.
\end{remark}

\section{Sketch of the new compactness method}

\medskip

   We present in the section the tool which has been used in \cite{BrJa} and  which is  the cornerstone to prove compactness on the density. The interested reader is also referred to \cite{BeJa}, \cite{BoBrMi},
 \cite{Po} for more on the corresponding critical spaces. This tool is really appropriate to cover more general equation of state or stress tensor form
 compared to the more standard defect measure criterion used in \cite{Li}, \cite{Fe},
 \cite{FeNo}, \cite{NoSt} for instance. 
%
\subsection{The compactness criterion}
We start by a well known result providing compactness of a sequence
\begin{proposition} \label{propcomp}
Let $\rho_k$ be a sequence uniformly equi-integrable in some $L^p((0,T)\times {\mathds T}^d)$ with $1\leq p<\infty$. Assume that ${\cal K}_h$ is a sequence of smooth, positive, bounded functions s.t.
\begin{eqnarray}
&i.\quad
\displaystyle  \forall \eta>0,\quad \sup_h \int_{{\mathds T}^d } {\cal K}_h(x) \, 
    \mathds{1}_{\{x\, : \, |x|\geq \eta\}}\,dx<\infty,\\
&ii.\quad \|{\cal K}_h\|_{L^1({\mathds T}^d)}\longrightarrow +\infty \qquad  \hbox{ as } h \to 0.
\end{eqnarray} 
Assume that $\partial_t \rho_k \in L^q(0,T,W^{-1,q}({\mathds T}^d))$  (with $q>1$)  uniformly in $k$ and
\begin{equation}
\limsup_{k}\,\sup_{t\in [0,T]} \Bigl[ \frac{1}{\|{\cal K}_h\|_{L^1}}\, \int_{{\mathds T}^{2d}} {\cal K}_h(x-y) \,|\rho_k(t,x)-\rho_k(t,y)|^p\,dx\,dy\Bigr] \longrightarrow 0,\ \mbox{as}\ h\rightarrow 0,\label{limsup=0}
\end{equation}
then $\rho_k$ is compact in $L^p ((0,T)\times{\mathds T}^d)$. Conversely if $\rho_k$ is compact in $L^p ((0,T)\times {\mathds T}^d)$ then the above quantity converges to $0$ as $h$ goes to zero.
\label{kernelcompactness}
\end{proposition}
For reader's convenience, we just quickly recall why (\ref{limsup=0}) implies the compactness in space (by simply
forgetting the time dependency). 
Denote $\bar {\cal K}_h$ the normalized kernel
\[
\bar {\cal K}_h=\frac{{\cal K}_h}{\|{\cal K}_h\|_{L^1}}.
\] 
Write
\begin{eqnarray}
\|\rho_k - \bar {\cal K}_h\star_x \rho_k \|_{L^p}^p 
& &\displaystyle \le\frac{1}{ \|{\cal K}_h\|_{L^1}^p}\int_{{\mathds T}^d}\Bigl(\int_{{\mathds T}^d} 
    {\cal K}_h(x-y) |\rho_k(t,x)-\rho_k(t, y)|dx \Bigr)^p \, dy\\
&&\displaystyle \nonumber \le
  \frac{1}{ \|{\cal K}_h\|_{L^1}}\int_{{\mathds T}^{2d}} {\cal K}_h(x-y) |\rho_k(t,x)-\rho_k(t, y)|^p dx \, dy,
\end{eqnarray}
which converges to zero uniformly in $k$ as the limsup is $0$ for the sup in time. 
On the other-hand for a fixed $h$, $\overline {\cal K}_h\star_x u_k$ is compact in $k$ so for example for any $z>0$
\begin{eqnarray}
\|\rho_k(\cdot)  - \rho_k(\cdot+z) \|_{L^p}& &\leq 2\,\|\rho_k - \bar {\cal K}_h\star_x \rho_k \|_{L^p}+\|\bar {\cal K}_h\star_x\rho_k - \bar {\cal K}_h\star_x\rho_k(.+z)\|_{L^p}\\
& &\leq 2\,\|\rho_k - \bar {\cal K}_h\star_x \rho_k \|_{L^p}+|z|\,\|\rho_k\|_{L^p}\,\|\bar {\cal K}_h\|_{W^{1,\infty}}.
\end{eqnarray}
This shows by optimizing in $h$ that
\[
\sup_k\|\rho_k(\cdot) - \rho_k(.+z) \|_{L^p}\longrightarrow 0,\quad\mbox{as}\ |z|\rightarrow 0.
\]
proving the compactness in space by the Rellich criterion. 

Concerning the compactness in time, one only has to use the uniform bound on $\partial_t \rho_k$
in $L^q(0,T;W^{-1,q}({\mathds T}^d)$ with $q>1$ . Taking any convolution kernel $L_\eta$ this implies that $\|L_\eta\star_x\rho_k\|_{W^{1,q}_{t,x}}\leq C\,\eta^{-\theta}$ for some exponent $\theta$ (where we only convolve in space). Therefore for any fixed $\eta$,  $L_\eta\star_x\rho_k$ is compact in $L^q$ and in fact compact in $L^p$ thanks to the equi-integrability of $\rho_k$ in that space. Extracting a converging subsequence, one has that $\|L_\eta\star_x \rho_k-L_\eta\star_x \rho\|_{L^p_{t,x}}\rightarrow 0$
as $k$ goes to $+\infty$. Now simply write
\begin{eqnarray}
\|\rho_k-\rho\|_{L^p_{t,x}}&\leq& \|L_\eta\star_x \rho_k-L_\eta\star_x \rho\|_{L^p_{t,x}}+\|\rho_k-L_\eta\star_x\rho_k\|_{L^p_{t,x}}+\|\rho-L_\eta\star_x\rho\|_{L^p_{t,x}}\nonumber \\
&\leq &
\|L_\eta\star_x \rho_k-L_\eta\star_x \rho\|_{L^p_{t,x}} \nonumber \\
&& +2\,\sup_k\|\int L_\eta(z)(\rho_k(t, .) - \rho_k(t, .+z)) \, dz \|_{L^p_{t,x}}\longrightarrow 0,
\end{eqnarray}
by optimizing in $\eta$.

\bigskip

%

\noindent 
{\it The ${\cal K}_{h_0}$ functions.} Define $K_h$ a sequence of non negative functions, 
\[
K_h(x)=\frac{1}{(h+|x|)^a},\quad \mbox{for}\ |x|\leq 1/2,
\]
with some $a>d$ and $K_h$ non negative,
independent of $h$ for $|x|\geq 2/3$, with support in $B(0,3/4)$ and periodized such as to belong in
$C^\infty({\mathds T}^d\setminus B(0,3/4))$. 

For convenience, let us denote
$$
\overline K_h(x)=\frac{K_h(x)}{\|K_h\|_{L^1}}.
$$
    For $0<h_0<1$, the important quantity to be used in Proposition \ref{propcomp} will be
$$
{\cal K}_{h_0}(x)=\int_{h_0}^1 \overline K_h(x)\,\frac{dh}h
$$
where
$$
K_h(x) =\frac{1}{(h+|x|)^a},\quad \mbox{for}\ |x|\leq 1/2. 
$$
Remark  the important property: $\|{\cal K}_{h_0}\|_{L^1} \sim |\log h_0|$.

\section{Proof of  Theorem \ref{MainResultPressureLaw} }
  As usually the proof of global weak solutions of PDEs is divided in three steps:
\begin{itemize}
\item {\it A priori} energy estimates and control of unknowns,
\item Stability of weak sequences: Compactness,
\item Construction of approximate solutions.
\end{itemize}

\subsection{Energy estimates and control of unknowns.}

\noindent {\it Energy estimate.}
  Le us multiply the Stokes equation by $u$ and integrate by parts, we get
$$\mu \int_{{\mathds T}^d} |\nabla u_k|^2  + \alpha \int_{{\mathds T}^d} |u_k|^2
    + \int_{{\mathds T}^d} \nabla P(\rho_k) \cdot u = \int_{{\mathds T}^d} S_k\cdot u_k.$$
Now we write the equation satisfied by $\rho_k e(\rho_k)$ where 
$e(\rho_k) = \int_{\rho_{\rm ref}}^{\rho_k} P(s)/s^2 ds$, with $\rho_{\rm ref}$ a constant
reference density, we get
$$\partial_t(\rho_k e(\rho_k)) + {\rm div}(\rho_k e(\rho_k) u_k) + P(\rho_k) {\rm div} u_k = 0.$$
   Integrating in space  and adding to the first equation we get
$$\frac{d}{dt} \int_{{\mathds T}^d} \rho_k e(\rho_k) 
    +   \mu \int_{{\mathds T}^d} |\nabla u_k|^2 + \alpha \int_{{\mathds T}^d} |u_k|^2    = 
    \int_{{\mathds T}^d} S_k\cdot u_k.$$ 
This gives the following estimate 
\begin{eqnarray}\label{est}
&& \sup_{t\in [0,T]} \int_{{\mathds T}^d} [\rho_k e(\rho_k)] (t)
   +    \int_0^T \int_{{\mathds T}^d}
     (\mu |\nabla u_k|^2  + \alpha |u_k|^2) \\
&  & \nonumber \hskip4cm   
     =   \int_0^T \int_{{\mathds T}^d} S_k\cdot u_k
       +  \int_{{\mathds T}^d} (\rho_k)_0 e((\rho_k)_0) 
\end{eqnarray} 
Assuming $((\rho_k)_0)^\gamma \in L^\infty(0,T;L^1({\mathds T}^d))$ uniformly,  one only needs the right-hand side quantity   $S_k \in L^2([0,T], H^{-1}({\mathds T}^d))$ uniformly. Using the behavior of $P$
we get the uniform bound
$$\rho_k^\gamma \in L^\infty(0,T;L^1({\mathds T}^d)), \qquad
     u_k \in L^2(0,T;H^1({\mathds T}^d)).$$

\noindent {\it Remark.} Note that Relation (\ref{est}) is replaced, at the level of the global weak solutions, by the energy inequality
\begin{eqnarray}\label{ineg}
&& \sup_{t\in [0,T]} \int_{{\mathds T}^d} [\rho e(\rho)] (t)
   +    \int_0^T \int_{{\mathds T}^d}
     (\mu |\nabla u|^2  + \alpha |u|^2) \\
&  & \nonumber \hskip4cm   
\le   \int_0^T \int_{{\mathds T}^d} S\cdot u
       +  \int_{{\mathds T}^d} \rho_0 e(\rho_0) 
\end{eqnarray}

\bigskip

\noindent {\it Extra integrability on $\rho_k$.}
 When now considering the compressible system (\ref{Biononmono}), the divergence ${\rm div} u_k$
 is given 
 $${\rm div} u_k =  \frac{1}{\mu} P(\rho_k) + \frac{1}{\mu}\Delta^{-1} {\rm div} R_k$$ 
 with $R_k = S_k - \alpha u_k$.
    Therefore, since $\rho_k \in L^\infty(0,T;L^\gamma({\mathds T}^d))$, if we multiply by $\rho_k^\theta$, we get
 $$ I =   \int_0^T \int_{{\mathds T}^d} P(\rho_k) \rho_k^\theta
         = \mu  \int_0^T\int_{{\mathds T}^d} {\rm div} u_k \rho_k^\theta  
          - \int_0^T\int_{{\mathds T}^d} \Delta^{-1} {\rm div} R_k \,  \rho_k^\theta       
$$
which is easily bounded as follows
$$ I
         \le  \bigl[\mu \|{\rm div} u_k\|_{L^2((0,T)\times{{\mathds T}^d})} 
           + \|\Delta^{-1}{\rm div} R_k\|_{L^2((0,T)\times{{\mathds T}^d})} \bigr]
                                            \|\rho_k^\theta\|_{L^2((0,T)\times {{\mathds T}^d})}
 $$
 Thus using the behavior of $P$ and information on $u_k$ and $R_k$, we get for large density
 $$ \int_0^T \int_{{\mathds T}^d} (\rho^{\gamma+\theta})  \le 
         C + \varepsilon  \int_0^T \int_{{\mathds T}^d} (\rho^{2\theta}).$$
Thus we get a control on $\rho_k^{\gamma+\theta}$ if $\theta\le \gamma$.
Therefore, we get $\rho_k \in L^p((0,T)\times {{\mathds T}^d})$ with $p>2$ is 
$\gamma>1$.

\begin{remark} Note that for the barotropic compressible Navier-Stokes equations, we get
$$ \frac{1}{2}\frac{d}{dt} \int_{{\mathds T}^d} \rho |u_k|^2 + \frac{d}{dt} \int_{{\mathds T}^d} \rho_k e(\rho_k) 
    +   \mu \int_{{\mathds T}^d} |\nabla u_k|^2   = 0.$$
and
$$\int_0^T \int_{{\mathds T}^d} \rho_k^{\gamma+\theta} <+ \infty
$$
for $\theta \le 2\gamma/d - 1$ where $d$ is the space dimension. The constraint on $\gamma$ 
in \cite{BrJa} is different because of more restrictive integrability information (due to the presence of the total time derivative).
\end{remark}

\subsection{Stability of weak sequences: Compactness}
  We will prove the following result which is the main part of the proof
\begin{proposition}
Assume $(\rho_k,u_k)$ satisfy system  {\rm (\ref{Biononmono})} in a weak sense with a 
pressure law satisfying {\rm (\ref{gammacontrolBrinkman})}--{\rm (\ref{nonmonotonehypBrinkman})}
and with the following weak regularity
$$
\sup_k \|\rho_k^\gamma\|_{L^\infty_t\,L^1_x}<\infty,\qquad \sup_k \|\rho_k\|_{L^p_{t,x}}<\infty\quad\mbox{with}\ p\leq 2\gamma,
$$ 
and
$$\sup_k \|u_k\|_{L^2_t H^1_x}<\infty.$$
  If the source term $S_k$ is compact in $L^2([0,\ T],\ H^{-1} ({\mathds T}^d))$ and the initial density sequence $(\rho_k)_0$ is assumed to be compact and hence satisfies
\[
\limsup_{k} \Bigl[\frac{1}{\|K_h\|_{L^1}} \int_{{\mathds T}^{2d}} K_h(x-y) \bigl|(\rho_k^x)_0 - (\rho_k^y)_0\bigr| \Bigr] 
 = \epsilon(h) \to 0  \hbox{ as } h\to 0,
\]
then $\rho_k$ is compact in $L^q((0,T)\times {\mathds T}^d)$ for all $q<p$. \label{proproughsketch}
\end{proposition}
\begin{remark}
Here and in the following, we use the convenient notation 
$(\rho_k^x,u_k^x) =(\rho_k(t,x), u_k(t,x))$, $(\rho_k^y,u_k^y)=(\rho_k(t,y),u_k(t,y))$ 
and $(\rho_k^x)_0=\rho_k(t=0,x)$, $(\rho_k^y)_0=\rho_k(t=0,y)$.
\end{remark}

\begin{proof}  As mentioned in \cite{BrJa},  regulartity estimates of $\rho_k$ solution of a transport equation have been derived by G. Crippa and C. De Lellis in \cite{CD} using explicit control
on characteristics. 
   But we know that due to the weak regularity of ${\rm div} u_k$ (due to the coupling 
between  ${\rm div} u_k$ and $\rho_k$), we cannot expect to simply propagate the regularity assumed on the density. The idea is to accept to lose some of it by introducing appropriate weights.
and by working at the PDE level instead of the ODE level. More precisely, we consider
weights $w_k$ such that $w_k|_{t=0} =1$ and thus in particular, since $\rho_k^0$ is compact
 $$
\limsup_{k} \Bigl[\frac{1}{|\log h_0|} \int_{{\mathds T}^{2d}} {\cal K}_{h_0} (x-y)
   \bigl|(\rho_k^x)_0 - (\rho_k^y)_0\bigr| 
  ((w_k^x)_0 + (w_k^y)_0) \Bigr]  \to 0 \hbox{ as } h_0\to 0.
$$
Remark that 
$$\frac{1}{|\log h_0|} \int_{h_0}^1 \frac{\epsilon(h)}{h} dh \to 0 \hbox{ when } h_0 \to 0.$$
Let us now choose a weight satisfying a PDE which is dual to the continuity equation
\begin{equation}\label{weight}
\left\{
\begin{array}{rl}
& \partial_t w_k + u_k \cdot \nabla w_k + \lambda D_k w_k = 0, \\
& w_k\vert_{t=0} = (w_k)_0= 1,
\end{array}
\right.
\end{equation} 
with $\lambda$ a constant parameter to be chosen later on and an appropriate
positive damping terms $D_k$ which will depend on the unknowns
$(\rho_k,u_k)$ and chosen also later-on.
We denote as before $w_k^x = w_k(t,x)$ and $w_k^y = w_k(t,y)$. It is convenient for the calculation to write the two equations for $w_k^x$ and $w_k^y$ (even though formally this is only Eq. (\ref{weight}))
\begin{equation}\label{weightx}
\left\{
\begin{array}{rl}
& \partial_t w_k^x + u_k ^x\cdot \nabla_x w_k^x + \lambda D_k^x w_k^x = 0, \\
& w_k^x\vert_{t=0} = (w_k^x)_0= 1,
\end{array}
\right.
\end{equation} 
and 
\begin{equation}\label{weighty}
\left\{
\begin{array}{rl}
&\partial_t w_k^y + u_k ^y\cdot \nabla_y w_k^y + \lambda D_y^x w_k^y = 0, \\
& w_k^y\vert_{t=0} = (w_k^y)_0=  1.
\end{array}
\right.
\end{equation}
 We first  study the propagation of the quantity 
$$
R_{h_0} (t)= \int_{{\mathds T}^{2d}} {\cal K}_{h_0}(x-y)\,
                   \bigl|\rho_k^x - \rho_k^y\bigr|\,(w^x+w^y)\,dx\,dy
                  = \frac{1}{\|K_h\|_{L^1}} \int_{h_0}^1 R(t) \frac{dh}{h}
$$
where 
$$
R(t) =  \int_{{\mathds T}^{2d}} K_{h}(x-y)\,
    \bigl|\rho_k^x - \rho_k^y \bigr|\,(w^x+w^y)\,dx\,dy.
$$
We show that it is possible to choose $D_k$ and $\lambda$
such that
 $$
\limsup_{k} \Bigl[\frac{1}{|\log h_0|} \int_{{\mathds T}^{2d}} {\cal K}_{h_0} (x-y) 
     \bigl|\rho_k^x - \rho_k^y\bigr| 
  (w_k^x + w_k^y)\Bigr]   \to 0 \hbox{ as } h_0\to 0
$$
as initially. Then,  we will need properties on $w_k^x$ (and hence $w_k^y$) to conclude that
we also have 
 $$
\limsup_{k} \Bigl[\frac{1}{|\log h_0|} \int_{{\mathds T}^{2d}} {\cal K}_{h_0} (x-y)
   \bigl|\rho_k^x - \rho_k^y\bigr| \Bigr]  \to 0 \hbox{ as } h_0\to 0
$$
which is the criterion giving compactness. Thus the proof is divided in two parts.

\medskip

\noindent {I) \em First step: Propagation of a weighted regularity.} 
Using the transport equation, we obtain that
\begin{eqnarray}
&\nonumber \partial_t |\rho_k^x-\rho_k^y |+{\rm div}_x\,(u_k^x \,|\rho_k^x-\rho_k^y|)
+{\rm div}_y\,(u_k^y\,|\rho_k^x-\rho_k^y|)\\
&\quad = \frac{1}{2}({\rm div}_x u_k^x+{\rm div}_y u_k^y)\,|\rho_k^x-\rho_k^y| 
-  \frac{1}{2}({\rm div}_x u_k^x- {\rm div}_y u_k^y)\,(\rho_k^x+\rho_k^y)\,s_k,
\end{eqnarray}
where $s_k=sign\,(\rho_k^x-\rho_k^y)$. Remark that these calculations can be justified
for a fixed $k$ through the DiPerna-Lions theory on renormalized solutions because
the densities and the gradient of the velocity are in $L^2$ in space and time.
    From this equation on $|\rho_k^x-\rho_k^y|$, we deduce by symmetry that 
\begin{eqnarray}
\frac{d}{dt}R(t) = &&\int_{{\mathds T}^{2d}} \nabla K_h(x-y)\cdot (u_k^x-u_k^y )\,|\rho_k^x-\rho_k^y|\,
     (w^x+ \,w^y)\\
&&\nonumber  
    -  \int_{{\mathds T}^{2d}}  K_h(x-y)\,({\rm div} u_k^x-{\rm div} u_k^y)\,(\rho_k^x+\rho_k^y + (\rho_k^x-\rho_k^y)\,)\,s_k\, w^x\\
&&\nonumber  + 2 \int_{{\mathds T}^{2d}}  K_h(x-y)\,|\rho_k^x -\rho_k^y |\,\left(\partial_t w_k^x
       +u_k^x \cdot\nabla_x w^x+  {\rm div}_x u_k^x\,w_k^x\right) \\
= & & \nonumber  A_1 + A_2 + A_3.
\end{eqnarray}

\medskip

\noindent {\it First term.} The first term will lead to non symmetric contributions. By definition of $K_h$, we have
$$|z| |\nabla K_h(z)| \le C K_h(z).$$
We hence write
\begin{eqnarray}
A_1 = &&\int_{{\mathds T}^{2d}} \nabla K_h(x-y)\cdot (u_k^x-u_k^y)\,|\rho_k^y-\rho_k^y|\,(w_k^x+w_k^y)\\
&&\nonumber 
\leq C\,\int_{{\mathds T}^{2d}} K_h(x-y)\,(D_{|x-y|} u_k^x + D_{|x-y|} u_k^y)\,|\rho_k^x-\rho_k^y | w_k^x,
\end{eqnarray}
where we have used here the inequality
\[
|u(x)-u(y)|\leq C\,|x-y|\,(D_{|x-y|} u_k^x + D_{|x-y|} u_k^y),
\]
with
\[
D_{h} u_k^x =\frac{1}{h}\,\int_{|z|\leq h} \frac{|\nabla u_k^{x+z}|}{|z|^{d-1}}\,dz.
\]
 This inequality is fully described in Lemma \ref{diffulemma} in the appendix
 with a proof given in~\cite{BrJa}.
   The key problem is the $(D_h u_k^y) w_k^x$ term which one will have to control by the term
$M|\nabla u_k^x|w_k^x$ in the penalization. This is where integration over $h$ and the use of translation
properties of operator will be used. For that we will add and subtract an appropriate quantity to obtain a symmetric expression.
  
Denoting $z=x-y$, using Cauchy-Shwartz inequality and the uniform $L^2$ bound on $\rho_k$, we have 
\begin{eqnarray}
\int_{h_0}^1 \int_0^t \frac{A_1}{\|K_h\|_L^1} \frac{dh}{h}  
&&
\leq C \int_{h_0}^1 \int_0^t \int_{{\mathds T}^d}
        \overline {K_h}(z) \| D_{|z|} u_k(\cdot) -D_{|z|} u_k(\cdot + z) \|_{L^2} \frac{dh}{h} 
\\
&&\nonumber 
 + C   \int_0^t \int_{{\mathds T}^{2d}}
        {\cal K}_{h_0}(x-y) D_{|x-y|} u_k(x) \,|\rho_k^x- \rho_k^y|  \, w_k^x. 
\end{eqnarray}
Using Lemma \ref{compareDhmax} which bounds $D_{|x-y|} u_k^x$ by the Maximal operator $M\,|\nabla u_k|(x)$, we deduce that 
\begin{eqnarray}
\int_{h_0}^1\int_0^t \frac{A_1}{\|K_h\|_L^1} \frac{dh}{h}  
&&
\leq C \int_{h_0}^1 \int_0^t \int_{{\mathds T}^d}
        \overline {K_h}(z) \| D_{|z|} u_k(\cdot) -D_{|z|} u_k(\cdot + z) \|_{L^2} \frac{dh}{h} 
\\
&&\nonumber 
 + C   \int_0^t \int_{{\mathds T}^{2d}}
        {\cal K}_{h_0}(x-y) M |\nabla u_k^x| \,|\rho_k^x- \rho_k^y|  w_k^x. 
\end{eqnarray}
The second term will be absorbed
using the weight definition. But the first quantity has to be controlled using the property
of the translation of operator $D_h$ and for this reason, this calculation is critical as it is the one which imposes the scales in ${\cal K}_{h_0}$.

\medskip

\noindent {\it Second term.}  Use the relation between ${\rm div}\, u_k^x$ 
(respectively  ${\rm div}\, u_k^y$) with $\rho_k^x$ (respectively $\rho_k^y$), to obtain
\[
  A_2 =  - \frac{2}{\mu}   \int_{{\mathds T}^{2d}}  K_h(x-y)\,(P(\rho_k^x) - P(\rho_k^y))\,\rho_k^x\,s_k\, w^x
              + Q_h(t) \\
\]
where $Q_h(t)$ encodes the compactness in space of $\Delta^{-1} {\rm div}\, R_k$ and therefore
may be forgotten for simplicity as 
$$\frac{1}{|\log h_0|} \int_0^t \int_{{\mathds T}^{2d}} {\cal K}_{h_0} (x-y) Q_h(t) \to 0 \hbox{ as } h_0 \to 0,$$
as $R_k$ is compact in $L^2_t H^{-1}_x$ and hence $\Delta^{-1} {\rm div}\, R_k$ is compact in $L^2_{t,x}$ by the gain of one derivative. 
   
The bad term $P(\rho_k^y) w_k^x$ cannot a priori be bounded directly
with weights.  Hence we have to work a little on the expression $A_2$.
Recall first that the weight is positive because of min principle.

\smallskip

Let us remind that $w\geq 0$ by the maximum principle.

\medskip

\noindent -- Case 1: The case where $(P(\rho_k^x) - P(\rho_k^y)) (\rho_k^x- \rho_k^y) \ge 0$ and hence $(P(\rho_k^x) - P(\rho_k^y))\,s_k \ge 0$.
Then we have the right sign for the contribution namely a negative sign.

\smallskip

\noindent -- Case 2: The case  $(P(\rho_k^x) - P(\rho_k^y)) (\rho_k^x- \rho_k^y) < 0$ and
$\rho_k^y \le \rho_k^x/2$ or $\rho_k^y \ge 2 \rho_k^x$. 

\noindent a) Assume we are in the case $\rho_k^y \ge 2 \rho_k^x$ then, due to 
the assumption on the sign of the product $(P(\rho_k^x) - P(\rho_k^y)) (\rho_k^x- \rho_k^y))$, we have $P(\rho_k^x) > P(\rho_k^y)$. No we remark that
\[P(\rho_k^x) - P(\rho_k^y) = (P(\rho_k^x) + C )- (P(\rho_k^y)+C) \le P(\rho_k^x) + C\]
where $C$ is the constant in the lower bound of the pressure. Thus
\[P(\rho_k^x) - P(\rho_k^y)  \le C((\rho_k^x)^\gamma + 1)\]
since $P(\xi) \le P(0) + C \xi^{\gamma-1} \xi \le C\xi^\gamma.$ 
 Now remark that
$$ \rho_k^x \le \rho_k^y - \rho_k^x$$
and thus using the bound on $P(\rho_k^x)- P(\rho_k^y)$, we get
$$[P(\rho_k^x)- P(\rho_k^y)]\rho_k^x \le C((\rho_k^x)^\gamma +1)
    (\rho_k^y - \rho_k^x).$$
and then due to the negative sign of $\rho_k^x - \rho_k^y$ 
\[
(P(\rho_k^x) - P(\rho_k^y)) \rho_k^x s_k
    \ge  - C\,((\rho_k^x)^\gamma+1)\, |\rho_k^x - \rho_k^y|.
\]

\noindent b) Assume we are in the case $\rho_k^y \le \rho_k^x/2$ then, due to the assumption on the sign of the product $(P(\rho_k^x) - P(\rho_k^y)) (\rho_k^x- \rho_k^y))$, we have
 $P(\rho_k^x) < P(\rho_k^y)$.  As previously
 \[P(\rho_k^y) - P(\rho_k^x) = (P(\rho_k^y) + C )- (P(\rho_k^x)+C) \le P(\rho_k^y) + C\]
 and thus 
 \[P(\rho_k^y) - P(\rho_k^x)  \le C((\rho_k^y)^\gamma + 1) \le C((\rho_k^x)^\gamma +1)\]
Now using the trick
\[ \rho_k^x  \le \rho_k^x + \rho_k^x -2\rho_k^y = 2 (\rho_k^x - \rho_k^y) \] 
we get
\[
(P(\rho_k^x) - P(\rho_k^y)) \rho_k^x  s_k
 \ge  -  C\, ((\rho_k^x)^\gamma +1) \, |\rho_k^x - \rho_k^y|.
\]

\smallskip

\noindent --- Case 3: The case where $P(\rho_k^x) - P(\rho_k^y)$ and $\rho_k^x-\rho_k^y$ have different  signs but $\rho_k^x/2 \le \rho_k^y \le 2\rho_k^x$. Then we use the Lipschitz bound on
$p$ to get
\[ |P(\rho_k^x) - P(\rho_k^y)| \le C( (\rho_k^x)^{\gamma-1} 
   + (\rho_k^y)^{\gamma-1})|\rho_k^x- \rho_k^y|
   \le C (\rho_k^x)^{\gamma-1}|\rho_k^x- \rho_k^y|
\]
and thus
\[
(P(\rho_k^x) - P(\rho_k^y)) \rho_k^x  s_k
 \geq -C\,   (\rho_k^x)^\gamma \, |\rho_k^x - \rho_k^y|.
\] 

\medskip

\noindent Therefore we get the following interesting bound:
\[
A_2 \le C \int K_h(x-y)\, (1+ (\rho_k^x)^\gamma))\, |\rho_k^x-\rho_k^y|\, w_k^x.
\]

\medskip

\noindent {\it Third term.} Using the equations satisfied by $w_k^x$ and $w_k^y$, we
have 
\begin{eqnarray}
A_3 = && \int_{{\mathds T}^{2d}}  K_h(x-y)\,|\rho_k^x -\rho_k^y |\,\left(\partial_t w_k^x
       +u_k^x \cdot\nabla_x w^x+ \,{\rm div}_x u_k^x\,w_k^x\right) \\
&&\nonumber  \hskip3cm      \le
\int_{{\mathds T}^{2d}}  K_h(x-y)\,|\rho_k^x -\rho_k^y |\,\left(-\lambda D_k^x  
   + \,{\rm div}_x u_k^x \right)w_k^x.
\end{eqnarray}

\medskip

\noindent {\it Conclusion of the first step.}
Collecting the three steps, we get
\begin{eqnarray}
R_{h_0} (t)  - R_{h_0}(0)&&\le C \int_{h_0}^1 \int_0^t \int_{{\mathds T}^{2d}}
        \overline {K_h}(z) \| D_{|z|} u_k(\cdot) -D_{|z|} u_k(\cdot + z) \|_{L^2} \frac{dh}{h}Ê\\
&&\nonumber   + C   \int_0^t \int_{{\mathds T}^{2d}}
        {\cal K}_{h_0}(x-y) M |\nabla u_k^x| \,|\rho_k^x- \rho_k^y|  \, w_k^x  \\
 &&  \nonumber  + C   \int_0^t \int_{{\mathds T}^{2d}}
       {\cal K}_{h_0}(x-y) (1+ (\rho_k^x)^\gamma)) |\rho_k^x-\rho_k^y| w_k^x \\
&& \nonumber +   \int_0^t \int_{{\mathds T}^{2d}}  {\cal K}_{h_0}(x-y)\,|\rho_k^x -\rho_k^y |\,\left(-\lambda D_k^x  + \,{\rm div}_x u_k^x \right)w_k^x.
\end{eqnarray}
Therefore we choose 
$$D_k = M|\nabla u_k| + |div u_k| + (\rho_k)^\gamma.$$
 Then for $\lambda$ large enough, we get
\begin{eqnarray}
R_{h_0} (t)  - R_{h_0}(0)&&\le C \int_{h_0}^1 \int_0^t \int_{{\mathds T}^{2d}}
        \overline {K_h}(z) \| D_{|z|} u_k(\cdot) -D_{|z|} u_k(\cdot + z) \|_{L^2} \frac{dh}{h}Ê\\
 &&  \nonumber  + C   \int_0^t R_{h_0}(\tau)  d\tau.
\end{eqnarray}
We now use translation property implied by the square functions given in 
Appendix, and more precisely using Lemma \ref{shiftDulemma} (proved in
\cite{BrJa}), we may write
\begin{eqnarray}
R_{h_0} (t)  - R_{h_0}(0)&&\le C\, |\log h_0|^{1/2}\, \int_0^t \|u(\tau,.)\|_{H^1_x}\,d\tau
   + C  \, \int_0^t R_{h_0}(\tau)  d\tau.
\end{eqnarray}
Therefore using that $u_k$ is uniformly bounded in $L^2(0,T;H^1({\mathds T}^d))$ and using the assumption on $R_{h_0}(0)$, then by Gronwall Lemma, we get that 
\[
\limsup_k \sup_{t\in[0,\ T]}\frac{R_{h_0}}{|\log h_0|}\longrightarrow 0,\quad\mbox{as}\ h_0\rightarrow 0, 
\] 
which is the desired propagation property.

\bigskip
      
\noindent{II) \em Second step.}  We now have to control the weights so as to
remove them. Namely we want to prove that
$$
\limsup_{k} [\frac{1}{|\log h_0|} \int_{{\mathds T}^{2d}} {\cal K}_{h_0}(x-y)\,|\rho_k^x - \rho_k^y | \,dx\,dy]
\to 0 \hbox{ as } h_0\to 0$$ 
and not only  
$$
\limsup_{k} [\frac{1}{|\log h_0|} \int_{{\mathds T}^{2d}} {\cal K}_{h_0}(x-y)\,|\rho_k^x - \rho_k^y | \,
     (w_k^x+w_k^y)dx\,dy] \to 0 \hbox{ as } h_0 \to 0.
$$
Remark that from its equation, the weight also satisfies 
$$\partial_t |\log w_k| + u_k \cdot \nabla |\log w_k| = \lambda D_k,$$
with
$$D_k = M|\nabla u_k| + |{\rm div} u_k| + (\rho_k)^\gamma.$$
Thus multiplying by $\rho_k$ and using the mass or continuity equation, we get
$$\frac{d}{dt} \int_{{\mathds T}^d} \rho |\log w_k| = \lambda  \int_{{\mathds T}^d} \rho D_k.$$
Note that $(u_k)_{k\in N}$ and $(\rho_k)_{k\in N}$ are respectively uniformly bounded in
 $L^2(0,T;H^1({\mathds T}^d))$ and $L^{2\gamma}((0,T)\times {\mathds T}^d)$ with
$\gamma > 1$, thus the right-hand side is uniformly bounded.

\bigskip

\noindent  Denoting $\omega = \{x: w_k \le \eta\}$, note that
\begin{eqnarray*}
\int_{{\mathds T}^{2d}} {\cal K}_{h_0}(x-y)\,|\rho_k^x - \rho_k^y | \,dx\,dy
&& = \int_{h_0}^1 \int_{{\mathds T}^{2d}} \overline K_h(x-y) |\rho_k^x - \rho_k^y| \frac{dh}{h}\\
&& = \int_{h_0}^1 
  \int_{x\in \omega_\eta^c \hbox{ or } y \in \omega_\eta^c} \overline K_h(x-y) |\rho_k^x - \rho_k^y| \frac{dh}{h} \\
&& + \int_{h_0}^1 
  \int_{x\in \omega_\eta \hbox{  and } y \in \omega_\eta} \overline K_h(x-y) |\rho_k^x - \rho_k^y| \frac{dh}{h}  \\
  && = B_1 + B_2.
\end{eqnarray*}
It suffices to observe that 
$$B_1 \le   \frac{1}{\eta} R_{h_0}$$
while 
by the property of the weights $w_k$ 
\[
B_2 \le 2 \int_{h_0}^1 \int_{{\mathds T}^{2d}} \overline K_h(x-y)\,\rho_k\, 1_{w_k\le \eta}\, \frac{dh}{h}\leq C\, \frac{|\log h_0|}{|\log \eta|}\,\int_{{\mathds T}^d} \rho_k\,|\log w_k|\,dx
             \le C\, \frac{|\log h_0|}{|\log \eta|}.
\]
 Combining the estimates, one obtains
 $$\int_{{\mathds T}^{2d}} {\cal K}_{h_0}(x-y)\,|\rho_k^x - \rho_k^y | \,dx\,dy
 \le C\Bigl( \frac{\int_{h_0}^1 \varepsilon(h) \frac{dh}{h} + |\log h_0|^{1/2}}{\eta}   
        + \frac{\|{\cal K}_{h_0}\|_{L^1}}{|\log\eta|} \Bigr) 
 $$
 and therefore
 \begin{eqnarray*}
& &\frac{1}{\|{\cal K}_{h_0}\|_{L^1}} \int_{{\mathds T}^{2d}} {\cal K}_{h_0}(x-y)\,|\rho_k^x - \rho_k^y | \,dx\,dy
\\
& &\qquad \le
 \displaystyle  C\left( \frac{ \displaystyle \frac{1}{|\log h_0|} \int_{h_0}^1 \varepsilon(h) \frac{dh}{h} 
       +  \displaystyle  |\log h_0|^{-1/2}}{\eta}   
        + \frac{1}{|\log\eta|} \right). 
 \end{eqnarray*}
 Denoting $\overline\epsilon(h_0) = \int_{h_0}^1 \epsilon(h)/h \, dh$ and optimizing $\eta$, we get 
  $$\frac{1}{\|{\cal K}_{h_0}\|_{L^1}} \int_{{\mathds T}^{2d}} {\cal K}_{h_0}(x-y)\,|\rho_k^x - \rho_k^y | \,dx\,dy
 \le
 \displaystyle   \frac{C}{|\log \bigl(|\log h_0|^{-1/2}   + \overline\epsilon(h_0)) \bigr)|^{1/2}}.
 $$ 
   This control in terms of $h_0$ coupled with  the uniform bound on $\partial_t \rho_k$ we get  
 using the mass equation and the estimates coming from the energy allows to apply the 
 compactness lemma and conclude that   $\rho_k$ is compact in $L^1((0,T)\times {\mathds T}^d)$.
   Thus $\rho_k$ is compact in  $L^q(0,T)\times {\mathds T}^d)$ for all $q<p$ using the extra integrability on
$\rho_k$. This gives the compactness property to pass to the limit in the non-linear terms.
 \bigskip
 
 \noindent {\it Remark.}
  The choice of appropriate weights is important in the proof.  It really depends on the
 system under consideration.  In \cite{BrJa}, we can find various choices depending on
 pressure laws or anisotropy in the viscous tensor. These weights penalize in some sense
 bad trajectories.

\end{proof}

\subsection{Construction of approximate solutions.}
  Our starting point for global existence is the following regularized system  
\begin{equation}\label{Barocompeps}
\left\{
\begin{array}{rl}
& \partial_t \rho_{k} + {\rm div} (\rho_{k} u_{k}) 
     = \alpha_k \Delta \rho_{k},\\
&      - \mu \Delta u_k - (\lambda+\mu) \nabla {\rm div } u_k 
  + \nabla P_\epsilon(\rho_k) + \alpha_k \nabla \rho_k\cdot \nabla u_k
                                    = S,
\end{array}
\right. 
\end{equation}
with the fixed source term $S$ and the fixed initial data
\begin{equation}
\rho_k\vert_{t=0} =\rho^0.\label{datainitial}
\end{equation} 
The pressure $P_\epsilon$ is defined as follows:
$$P_\epsilon(\rho) = p(\rho) \hbox{ if } \rho \le c_{0,\epsilon}, \qquad
    P_\epsilon (\rho)= p(C_{0,\epsilon}) + C(\rho - c_{0,\epsilon})^\beta
     \hbox{ if } \rho \ge c_{0,\epsilon},$$
 with large enough $\beta$.
  As usual the equation of continuity is regularized by means of an artificial viscosity term
and the momentum balance is replaced by a Faedo-Galerkin approximation to eventually reduce the problem on $X_n$, a finite-dimensional vector space of functions. 

   This approximate system can then be solved by a standard procedure: The velocity $u_k$ of the approximate momentum equation is looked as a fixed point of a suitable integral operator. Then 
given $u_k$, the approximate continuity equation is solved directly by means of the standard theory of linear parabolic equations. This methodology concerning the compressible Navier--Stokes equations is well explained and described in the reference books \cite{FeNo}, \cite{NoSt}. We omit the rest of this classical (but tedious) procedure and we assume that we have well posed and smooth solutions to (\ref{Barocompeps})--(\ref{datainitial}). 

     We now use the classical energy and extra bounds estimates detailed in the previous section. Note that they remain the same in spite of the added viscosity in the continuity equation. This is the reason in particular for the added term $\alpha_k \nabla \rho_k\cdot \nabla u_k$ in the momentum equation to keep the same energy balance. Let us summarize the {\it a priori} estimates that are obtained
$$
\sup_{k,\epsilon}\,\sup_t \int_{{\mathds T}^d} \rho_k^\gamma\,dx<\infty,\qquad
\sup_{k,\epsilon}\int_0^T\int_{{\mathds T}^d} |\nabla u_k|^2\,dx\,dt<\infty,
$$
and 
$$
\sup_{k,\epsilon}\int_0^T\int_{{\mathds T}^d} \rho_k^{p}(t,x)\,dx\,dt<\infty
$$ 
for all $p\le 2\gamma$.
      From those bounds it is straightforward to deduce that $\rho_k\,u_k$ belong to $L^q_{t,x}$ for some $q>1$, uniformly in $k$ and $\epsilon$. Therefore using the continuity equation bounds on
$\partial_t \rho_k$. 
    We have now to show the compactness of $\rho_k$ in $L^1$ and we can use the procedure
mentioned in \cite{Fe1} letting $\alpha_k$ goes to zero.   
     Then extracting converging subsequences, we can pass to the limit in every term (by classical approach) and obtain the existence of weak solutions to
\begin{equation}\label{Barocompepslim}
\left\{
\begin{array}{rl}
& \partial_t \rho + {\rm div} (\rho u) 
     = 0,\\
&
      - \mu \Delta u  + \alpha u   + \nabla P_\epsilon(\rho) = S.
\end{array}
\right. 
\end{equation} 
It remains then to pass to the limit with respect to $\epsilon$. This is
done using the stability procedure developed in the previous subsection
concerning compactness for general pressure laws.
%
\section{The compressible Navier-Stokes equations}
    We state in this section the main existence results that have been obtained in \cite{BrJa}.
    There exist several differences and complications compared to the global existence result we proved in this short
paper due in particular to the presence of the total time derivative.
This leads to more restrictions on the coefficient $\gamma$ in the pressure law. 
It could be interesting to try to extend our results with better gamma exponent using the renormalization procedure in \cite{Fe}
or with anisotropy in the stress tensor. 
\bigskip

\noindent {\bf I) The isotropic compressible Navier--Stokes equations with general pressure laws.}
Let us consider the isotropic compressible Navier--Stokes equations  
\begin{equation}\label{BaroPressureLaw}
\left\{
\begin{array}{rl}
& \partial_t \rho + {\rm div} (\rho u) =0,\\
& \partial_t (\rho u) + {\rm div}(\rho u\otimes u) - \mu \Delta u - (\lambda+\mu) \nabla {\rm div} u + \nabla P(\rho) =\rho f,
\end{array}
\right. 
\end{equation} 
 with ${2}\,\mu/d+\lambda>0$, a pressure law $P$ which is continuous on $[0,+\infty)$, $P$ locally Lipschitz on $(0,+\infty)$ with $P(0)=0$ such that there exists $C>0$ with
\begin{equation}\label{gammacontrol}
C^{-1} \rho^\gamma - C  \le P(\rho) \le C \rho^\gamma + C
\end{equation}
and for all $s\ge 0$
\begin{equation}
|P'(s)|\leq \bar P s^{\tilde \gamma-1}.
\label{nonmonotonehyp}
\end{equation}
\medskip

One then has global existence
\begin{theorem} \label{MainResultPressureLaw1} 
Assume that the initial data $u_0$ and $\rho_0\ge 0$ with $\int_{{\mathds T}^d} \rho_0 = M >0$ satisfy  the bound
$$
E_0= \int_{{\mathds T}^d} \bigl(\rho_0\frac{|u_0|^2}{2} + \rho_0 e(\rho_0)\bigr)\,dx <+\infty.
$$
   Let the pressure law $P$  satisfies
{\rm (\ref{gammacontrol})} and {\rm (\ref{nonmonotonehyp})} with 
\begin{equation}
\gamma>\, \bigl(\max(2,\tilde \gamma) +1\bigr)\, \frac{d}{d+2}.\label{hypgammaiso}
\end{equation} 
   Then there exists a global  weak solution of the compressible Navier--Stokes system~{\rm (\ref{BaroPressureLaw})} with positive density satisfying the initial data conditions in 
   ${\cal D}'({\mathds T}^d)$:
 $$\rho\vert_{t=0} = \rho_0, \qquad \rho u\vert_{t=0} = \rho_0u_0.$$ 

\medskip

\noindent Moreover the solution satisfies the explicit regularity estimate 
$$
\int_{{\mathds T}^{2d}} 1_{\rho_k(x)\geq \eta}\,1_{\rho_k(y)\geq \eta}\,K_h(x-y)\,\chi(\delta\rho_k)\leq \frac{C\,\|K_h\|_{L^1}}{\eta^{1/2}\,
|\log h|^{\theta/2}},
$$
for some $\theta>0$ where $\chi$ is a $C^2$ function such that
$\chi(\xi) = |\xi|^2$ if $|\xi|\le 1/2$ and $\chi(\xi)=|\xi|$ if $|\xi|>1$.
\end{theorem}

\bigskip

\noindent {\bf II) A non-isotropic compressible Navier--Stokes equations.}
We consider an example of non-isotropic compressible Navier--Stokes equations
\begin{equation}\label{BaroAnisotropic}
\left\{
\begin{array}{rl}
& \partial_t \rho + {\rm div} (\rho u) =0, \\
& \partial_t (\rho u) + {\rm div}(\rho u\otimes u) - {\rm div} \,(A(t)\, \nabla u)  
- (\mu+\lambda) \nabla {\rm div} u + \nabla P(\rho) =0,
\end{array}
\right. 
\end{equation} 
 with $A(t)$ a given smooth and symmetric matrix, satisfying
\begin{equation}
A(t)=\mu\,Id+\delta A(t),\quad \mu>0,\quad \frac{2}d\,\mu+\lambda-\|\delta A(t)\|_{L^\infty}>0.
\end{equation}
 We again take $P$ continuous on $[0,+\infty)$ with $P(0)=0$  but require it to be monotone 
 after a certain point
\begin{equation}
C^{-1}\,\rho^{\gamma-1}-C\leq P'(\rho)\leq C\,\rho^{\gamma-1}+C.\label{strictmonotone}
\end{equation}
with $\gamma>d/2$.
The second main result that we obtain is 
\begin{theorem} \label{MainResultAniso}
Assume that the initial data $u_0$ and $\rho_0\ge 0$ with $\int_{{\mathds T}^d} \rho_0 = M >0$ satisfy  the bound
$$
E_0= \int_{{\mathds T}^d} \bigl(\rho_0\frac{|u_0|^2}{2} + \rho_0 e(\rho_0)\bigr)\,dx <+\infty.
$$
Let  the pressure $P$ satisfies (\ref{strictmonotone}) with 
$$
\gamma > \frac{d}2 \left[\left(1+\frac{1}d\right) + \sqrt{1+\frac1{d^2}}\right]. \label{hypgammanoniso}
$$
  There exists a universal constant $C_\star$
 such that if 
$$
\|\delta A\|_{\infty} \le C_\star\,(2\mu+\lambda),
$$
  then there exists a global weak solution of the compressible  Navier--Stokes equation  
{\rm (\ref{BaroAnisotropic})} with positive density satisfying the initial data conditions
in ${\cal D}'({\mathds T}^d)$:
 $$\rho\vert_{t=0} = \rho_0, \qquad \rho u\vert_{t=0} = \rho_0u_0.$$ 
 The isotropic energy inequality is replaced by the  following anisotropic energy
$$
E(\rho,u)(\tau) +\int_0^\tau \int_\Omega (\nabla_x u^T\,A(t)\,\nabla u
    + (\mu+\lambda)\, |{\rm div} u|^2) \le E_0.
$$
\end{theorem}
\section{Appendix}
  In this appendix, let us give different results which are used in the paper. 
  The interested reader  is referred to \cite{BrJa} for details and proofs but also
 \cite{Stein}.
  These concern Maximal functions, Square functions and translation of operators.    
   First we remind the well known inequality
\begin{equation}
|\Phi(x)-\Phi(y)|\leq C\,|x-y|\,(M|\nabla\Phi|(x)+M|\nabla \Phi|(y)),\label{maximalineq}
\end{equation}
where $M$ is the localized maximal operator
\begin{equation}
M\,f(x)=\sup_{r\leq 1} \frac{1}{|B(0,r)|}\,\int_{B(0,r)} f(x+z)\,dz.\label{defmaximal}
\end{equation}
 Let us mention several mathematical properties that may be proved, see \cite{BrJa}.
  First one has 
\begin{lemma} There exists $C>0$ s.t. for any $u\in W^{1,1}({\mathds T}^d)$, one has
\[
|u(x)-u(y)|\leq C\,|x-y|\,(D_{|x-y|}u (x)+D_{|x-y|}u (y)), 
\]
where we denote
\[
D_h u(x)=\frac{1}{h}\,\int_{|z|\leq h} \frac{|\nabla u(x+z)|}{|z|^{d-1}}\,dz.
\]\label{diffulemma}
\end{lemma}
   Note that this result implies the estimate (\ref{maximalineq}) as
\begin{lemma}
\label{compareDhmax}
 There exists $C>0$, for any $u\in W^{1,p}({\mathds T}^d)$ with $p\geq 1$
\[
D_h\,u(x)\leq C\,M|\nabla u|(x).
\]\
\end{lemma}

\bigskip

    The key improvement in using $D_h$ is that small translations of the operator $D_h$ are actually easy to control
\begin{lemma} \label{shiftDulemma}
  Let $u\in H^1({\mathds T}^d)$ then  have the following estimates
\begin{equation}
\int_{h_0}^{1} \int_{{\mathds T}^d} \overline K_h(z)\,\|D_{|z|}\,u(.)-D_{|z|}\,u(.+z)\|_{L^2}\,dz\,\frac{dh}{h}\leq C\,|\log h_0|^{1/2}\,\|u\|_{H^1}.\label{squarefunction}
\end{equation}
\end{lemma} 
This lemma is critical and explain why we propagate a quantity integrated 
with respect to $h$ with a weight $dh/h$ namely with the Kernel ${\cal K}_{h_0}$. The full proof is rather classical and can be found in \cite{BrJa} for any $L^p$ space but we give a brief sketch here (which is simpler as Lemma 
\ref{shiftDulemma} is $L^2$ based and we can use Fourier transform).

\begin{proof}[Sketch of the proof of Lemma \ref{shiftDulemma}.]
Note that we can write
\[
D_h\,u(x)=L_h\star \nabla u,\quad L(x)=\frac{ 1_{|x|\leq1}}{|x|^{d-1}},\quad L_h(z)=h^{-d}\,L(z/h),
\]
where $L_h$ is hence a usual convolution operator and $L\in W^{s,1}$ for any $s<1$. Now 
\begin{eqnarray}
&&\displaystyle\int_{h_0}^{1} \int_{{\mathds T}^d} \overline K_h(z)\,\|D_{|z|}\,u(.)-D_{|z|}\,u(.+z)\|_{L^2}\,dz\,\frac{dh}{h}\nonumber\\
&&\qquad\leq C\, \int_{S^{d-1}}\int_{0}^{1}  \|L_r\star \nabla u(.)-L_r\star \nabla u(.+r\,\omega)\|_{L^2}\,\frac{dr}{r+h_0}\,d\omega\nonumber\\
&&\qquad\leq C\,|\log h_0|^{1/2}\,\left(\int_{S^{d-1}}\int_{0}^{1}  
  \|L_r\star \nabla u(.)-L_r\star \nabla u(.+r\,\omega)\|_{L^2}^2\,\frac{dr}{r+h_0}\,d\omega\right)^{1/2}.\nonumber
\end{eqnarray}
For any $\omega\in S^{d-1}$, define $L_r^\omega=L_r(.)-L_r(.+r\,\omega)$. Calculate
by Fourier transform
\[
\int_{0}^{1}  \|L_r\star \nabla u(.)-L_r\star \nabla u(.+r\,\omega)\|_{L^2}^2\,\frac{dr}{r+h_0}=\int_{0}^{1}\sum_{\xi\in{\mathds T}^d} |\hat L_r^\omega|^2(\xi)\,|\xi|^2 |\hat u|^2(\xi)\,\frac{dr}{r+h_0}.
\]
But $\hat L_r^\omega=(1-e^{ir\,\xi\cdot\omega})\,\hat L(r\,\xi)$ and furthermore $|\hat L(r\,\xi)|\leq C\,(1+|r\,\xi|)^{-s}$ for some $s>0$ since $L\in W^{s,1}$. Therefore
\[
\int_{0}^{1}  |\hat L_r^\omega|^2(\xi)\,\frac{dr}{r+h_0}\leq C,
\] 
for some constant $C$ independent of $\xi$, $\omega$ and $h_0$. This is of course the famous square function calculation and lets us conclude. 
\end{proof}


\end{document}